\documentclass[12pt]{amsart}
\usepackage{amssymb}
\usepackage{amsfonts}
\usepackage{latexsym}
\usepackage{amsmath}

\theoremstyle{plain}
\newtheorem{theorem}{Theorem}

\theoremstyle{definition}
\newtheorem{definition}{Definition}
\newtheorem{remark}{Remark}

\numberwithin{equation}{section}
\numberwithin{theorem}{section}
\numberwithin{definition}{section}
\numberwithin{lemma}{section}
\numberwithin{corollary}{section}
\numberwithin{proposition}{section}
\numberwithin{remark}{section}
\numberwithin{example}{section}
\numberwithin{table}{section}

\def\R{\mathbb{R}}
\def\Z{\mathbb{Z}}
\def\N{\mathbb{N}}
\def\ds{\displaystyle}

\allowdisplaybreaks

\begin{document}

\title[Solutions of the Variational Equation]{Solutions of the Variational Equation for an nth Order Boundary Value Problem with an Integral Boundary Condition}

\author[B. L. Jeffers]{Benjamin L. Jeffers}

\address{Trinity University, 1 Trinity Place, San Antonio, TX 78212}
\email{bjeffers@trinity.edu}

\author[J. W. Lyons]{Jeffrey W. Lyons}
\address{The Citadel, 171 Moultrie Street, Charleston, SC 29409}
\email{jlyons3\symbol{64}citadel.edu}

\begin{abstract}  In this paper, we discuss differentiation of solutions to the boundary value problem $y^{(n)} = f(x, y, y^{'}, y^{''}, \ldots, y^{(n-1)}), \; a<x<b,\; 	y^{(i)}(x_j) = y_{ij},\; 0\leq i \leq m_j, \; 1 \leq j \leq k-1$, and $y^{(i)}(x_k) + \int_c^d p y(x)\;dx = y_{ik}, \;0 \leq i \leq m_k,\;\sum_{i=1}^km_i=n$ with respect to the boundary data. We show that under certain conditions, partial derivatives of the solution $y(x)$ of the boundary value problem with respect to the various boundary data exist and solve the associated variational equation along $y(x)$.
\end{abstract}

\maketitle\thispagestyle{empty}

{\bf MCS 2020:} 34B10, 34B15
 
\section{Introduction}
Our concern is characterizing partial derivatives with respect to the boundary data of solutions to the $n$th order nonlocal boundary value problem 
\begin{equation}\label{eq1}
	y^{(n)} = f\left(x, y, y^{'}, y^{''}, \ldots, y^{(n-1)}\right), \; a<x<b
\end{equation}
satisfying

\begin{equation}\label{eq2}
	\begin{array}{c}
		y^{(i)}\left(x_j\right) = y_{ij},\; 0\leq i \leq m_j, \; 1 \leq j \leq k-1, \\
		y^{(i)}\left(x_k\right) + \ds\int_c^d p y(x)\;dx = y_{ik}, \;0 \leq i \leq m_k
	\end{array}
\end{equation}
where and throughout $k,n\in\N$ with $2 \leq k \leq n,\; m_1, \ldots, m_k\in\Z^+$ such that $\sum_{i=1}^km_i=n,$ and $a < x_1 < x_2 < \cdots < x_k <c<d <b,p\in\R$.

Differentiation of solutions of initial value problems with respect to initial conditions has been a well-known result in the field of differential equations for a long time. In his book \cite{Hartman64}, Hartman attributes the theorem and proof to Peano. Hence, the result is commonly referred to as a theorem of Peano. These derivatives solve the associated variational equation to the differential equation.

Subsequently, similar results were obtained for boundary value problems and relied heavily upon the continuous dependence of solutions of boundary value problems on boundary conditions. The continuous dependence result utilizes a map of initial conditions to boundary conditions and the Brouwer Invariance of Domain Theorem. Results for boundary value problems on differential equations with standard boundary conditions may be found in \cite{Henderson84, Peterson76, Peterson78, Spencer75, Sukup75}. 

Direct analogues also exist for difference equations \cite{HendersonLee91} and dynamic equations on times scales \cite{BaxterLyonsNeugebauer16}. The mathematics community has added a parameter to the nonlinearity \cite{HendersonHornHoward94, HendersonJiang15, LyonsMiller15}. Researchers have also produced results for various types of boundary conditions including nonlocal \cite{BenchohraHamaniHendersonNtouyas07, EhrkeHendersonKunkelSheng07, HendersonTisdell04, HendersonHopkinsKimLyons08, HopkinsKimLyonsSpeer09, Lawrence02, Lyons11, Lyons14}, functional \cite{Datta98, Ehme93, EhmeEloeHenderson93, EhmeHenderson96, EhmeLawrence00}, and integral \cite{BenchohraHendersonLucaOuahab14,LyonsMajorSeabrook18}.

In this paper, we extend the results of \cite{JansonJumanLyons14} to an $n$th order differential equation using the procedure outlined in \cite{Henderson87}. The general idea is to use continuous dependence to write the solution of the boundary value problem as the solution to an initial value problem. After multiple applications of the Mean Value Theorem, we can apply Peano's theorem directly to the problem at hand.

The remainder of this paper is organized as follows. In section two, we present the boundary value problem and define its associated variational equation. We also introduce five hypotheses that are imposed upon the differential equation along with Peano's Theorem and the continuous dependence result. Our boundary value problem with integral condition analogue is found in section three.

\section{Assumptions and Background Theorems}

We establish a few conditions that are imposed upon (\ref{eq1}):
\begin{enumerate}
	\item[(i)] $f\left(x,y_1, \ldots,y_n\right): (a,b) \times \R^n \to \R$ is continuous,
	\item[(ii)] $\frac{\partial f}{\partial y_i} \left(x,y_1, \ldots,y_n\right): (a,b) \times \R^n \to\R$ is continuous, $i = 1,\ldots,n,$
	\item[(iii)] solutions of initial value problems for (\ref{eq1}) extend to $(a,b).$
\end{enumerate}

\begin{remark} Note that \textup{(iii)} is not a necessary condition but lets us avoid continually making statements about maximal intervals of existence inside $(a,b)$.  
\end{remark}

\noindent Next, the results discussed rely upon the definition of the variational equation which we present here.
\begin{definition}
Given a solution $y(x)$ of \textup{(\ref{eq1})} and for $i=1,2,\ldots,n,$ we define the \textit{variational equation along $y(x)$} by
\begin{equation}\label{var}
z^{(n)} = \sum_{i=1}^n\frac{\partial f}{\partial y_i} \left(x,y,y',\ldots,y^{(n-1)}\right) z^{(i-1)}.
\end{equation}
\end{definition}

Our aim is an analogue of the following theorem that Hartman \cite{Hartman64} attributes to Peano for (\ref{eq1}), (\ref{eq2}).

\begin{theorem}\label{peano}[A Peano Theorem]  Assume that, with respect to \textup{(\ref{eq1})},
conditions \textup{(i)-(iii)} are satisfied.  Let $x_0 \in (a,b)$ and $$y(x)
:= y\left(x, x_0, c_0 ,c_1,\ldots,c_{n-1}\right)$$ denote the solution of \textup{(\ref{eq1})}
satisfying the initial conditions $y^{(i)}\left(x_0\right) = c_i,\; 0\leq i\leq n-1.$
Then,

\begin{enumerate}
\item[(a)] for each $0\leq j\leq n-1$, $\alpha_j(x) := \frac{\partial
y}{\partial c_j}(x)$ exists on $(a,b)$ and is the solution of the variational equation \textup{(\ref{var})} along $y(x)$ satisfying the initial conditions
$$\alpha_j^{(i)}\left(x_0\right)  =  \delta_{ij}, \; 0\leq i \leq n-1.$$
\item[(b)] $\beta(x)
:= \frac{\partial y}{\partial x_0}(x)$ exists on $(a,b)$ and is the solution of the variational equation \textup{(\ref{var})} along $y(x)$ satisfying the initial conditions
$$\beta^{(i)}\left(x_0\right)  =  -y^{(i)}\left(x_0\right), \; 0\leq i \leq n-1.$$
\item[(c)] $\frac{\partial y}{\partial x_0} (x) = -\sum_{i=0}^{n-1} y^{(i)}(x_0) \frac{\partial
y}{\partial c_i} (x).$
\end{enumerate}
\end{theorem}

The next condition guarantees uniqueness of solutions of (\ref{eq1}), (\ref{eq2}) and is a nonlocal analogue of $(m_1,\ldots,m_k)$-disconjugacy. 

\begin{enumerate}
\item[(iv)] If, for $0\leq i\leq m_j-1, \; 1\leq j\leq k-1,$
$$y^{(i)}\left(x_j\right) = z^{(i)}\left(x_j\right),$$ 
and, for $0\leq i\leq m_k-1,$
$$y^{(i)}\left(x_k\right) + \int_c^d p y(x)\;dx = z^{(i)}\left(x_k\right) + \int_c^d pz(x)\;dx,$$
where $y(x)$ and $z(x)$ are solutions of (\ref{eq1}), then, on $(a,b),$ $$y(x) \equiv z(x).$$
\end{enumerate}

The last condition provides uniqueness of solutions of (\ref{var}) along all solutions of (\ref{eq1}) and again is a nonlocal analogue of $(m_1,\ldots,m_k)$-disconjugacy.

\begin{enumerate}
\item[(v)] Given a solution $y(x)$ of (\ref{eq1}), if, for $0\leq i\leq m_j-1, \; 1\leq j\leq k-1,$
$$u^{(i)}\left(x_j\right)=0,$$
and, for $0\leq i\leq m_k-1,$ 
$$u^{(i)}\left(x_k\right) + \int_c^d p u(x)\;dx = 0,$$
where $u(x)$ is a solution of (\ref{var}) along $y(x),$ then, on $(a,b)$, 
$$u(x) \equiv 0.$$
\end{enumerate}

We also make use of the following continuous dependence result for boundary value problems. A typical proof may be found in \cite{HendersonKarnaTisdell05}.

\begin{theorem}\label{contdep}[Continuous Dependence on Boundary Conditions] Assume \textup{(i)-(iv)} are satisfied with respect to
\textup{(\ref{eq1})}.  Let
$y(x)$ be a solution of \textup{(\ref{eq1})} on $(a,b)$. Then, there
exists a $\delta > 0$ such that, for 
$$\left|x_j - t_j\right| < \delta,\; 1\leq j\leq k,$$
$$\left|c-\xi\right|<\delta,\; \left|d-\Delta\right|<\delta,\;\left|p-\rho\right|<\delta,$$
$$\left|y^{(i)}\left(x_j\right) - y_{ij} \right| < \delta,\; 0\leq i\leq m_j-1, \; 1\leq j\leq k-1,$$
and
$$\left|y^{(i)}\left(x_k\right)+ \int_c^d p y(x)\;dx - y_{ik}\right| < \delta, \; 0\leq i\leq m_k-1,$$
there exists a unique solution $y_{\delta} (x)$ of \textup{(\ref{eq1})} such that
$$y^{(i)}_{\delta}\left(t_j\right) = y_{ij},\; 0\leq i\leq m_j-1, \; 1\leq j\leq k-1,$$
$$y_{\delta}^{(i)}\left(t_k\right) + \int_{\xi}^{\Delta} \rho y_{\delta}(x)\;dx = y_{ik},\; 0\leq i\leq m_k-1,$$
and, for $0\leq i\leq n-1,$ $\{y_\delta^{(i)}(x)\}$ converges uniformly to $y^{(i)}(x)$ as $\delta \to 0$ on $[\alpha,\beta]\subset(a,b)$.
\end{theorem}

\section{Analogue of Peano's Theorem}
\indent

In this section, we present our analogue to Theorem \ref{peano} stated in five parts.
\begin{theorem}\label{mainresult} Assume conditions \textup{(i)-(v)} are satisfied. Let $u(x)=\\u(x, x_1,\ldots,x_k,y_{01},\ldots,y_{m_k-1,k},p,c,d)$ be the solution of \textup{(\ref{eq1})} on $(a,b)$ satisfying
$$u^{(i)}\left(x_j\right)=y_{ij}, \; 0\leq i\leq m_j-1, \; 1\leq j\leq k-1,$$ and $$u^{(i)}\left(x_k\right)+\int_c^d pu(x)dx=y_{ik},\; 0\leq i\leq m_k-1.$$  Then,

\begin{enumerate}
	\item [$(a)$] for each $1\leq l\leq k-1$ and $0\leq r\leq m_l-1,\; Y_{rl}(x):=\frac{\partial u}{\partial y_{rl}}(x)$ exists on $(a,b)$ and is the solution of the variational equation
	\textup{(\ref{var})} along $u(x)$ satisfying the boundary conditions
	\begin{align*}
		&Y^{(i)}_{rl}\left(x_j\right) = 0, \; 0 \leq i \leq m_j -1,\; 1 \leq j \leq k-1,\; j \neq l \\
		&Y^{(i)}_{rl}\left(x_l\right) = 0,\;0 \leq i \leq m_l - 1,\; i \neq r \\
		&Y^{(r)}_{rl}\left(x_l\right) = 1 \\
		&Y^{(i)}_{rl}\left(x_k\right) + \int_c^dpY_{rl}(x)dx = 0,\; 0 \leq i \leq m_k - 1,
	\end{align*}
	
	and for $0\le r \le m_k-1,\; Y_{rk}:=\frac{\partial u}{\partial y_{rk}}(x)$ exists on $(a,b)$ and is the solution of the variational equation \textup{(\ref{var})} along $u(x)$ satisfying the boundary conditions
	\begin{align*}
		&Y^{(i)}_{rk}\left(x_j\right) = 0,\; 0 \leq i \leq m_j -1,\; 1 \leq j \leq k-1, \\
		&Y^{(i)}_{rk}\left(x_k\right) + \int_c^dpY_{rk}(x)dx = 0,\; 0 \leq i \leq m_k - 1, \; i\neq r,\\
		&Y^{(r)}_{rk}\left(x_k\right) + \int_c^dpY_{rk}(x)dx = 1,
	\end{align*}
	
	\item [(b)] for each $1\leq l\leq k-1,\; X_{l}(x):=\frac{\partial u}{\partial x_{l}}(x)$ exists on $(a,b)$ and is the solution of the variational equation
	\textup{(\ref{var})} along $u(x)$ satisfying the boundary conditions
	\begin{align*}
		&X^{(i)}_{l}\left(x_j\right) = 0,\;0\le i\le m_j-1,\; 1 \leq j \leq k-1,\; j \neq l \\
		&X^{(i)}_{l}\left(x_l\right) = -u^{(i+1)}(x_l),\; 0 \leq i \leq m_l - 1, \\
		&X^{(i)}_{l}\left(x_k\right) + \int_c^dpX_{l}(x)dx = 0,\; 0 \leq i \leq m_k - 1,
	\end{align*}
	
	and $X_k:=\frac{\partial u}{\partial x_k}(x)$ exists on $(a,b)$ and is the solution of the variational equation \textup{(\ref{var})} along $u(x)$ satisfying the boundary conditions
	\begin{align*}
		&X^{(i)}_{k}\left(x_j\right) = 0,\; 0\le i\le m_j-1,\;1 \leq j \leq k-1, \\
		&X^{(i)}_{k}\left(x_k\right) + \int_c^dpX_{k}(x)dx = -u^{(i+1)}\left(x_k\right),\; 0 \leq i \leq m_k - 1.
	\end{align*}
	
	\item [(c)] $C(x):=\frac{\partial u}{\partial c}(x)$ exists on $(a,b)$ and is the solution of the variational equation \textup{(\ref{var})} along $u(x)$ satisfying the boundary conditions
	\begin{align*}
		&C^{(i)}\left(x_j\right) = 0,\; 0\le i\le m_j-1,\;1 \leq j \leq k-1, \\
		&C^{(i)}\left(x_k\right) + \int_c^dpC(x)dx = -pu(c),\; 0 \leq i \leq m_k - 1.
	\end{align*}
	
	\item [(d)]  $D(x):=\frac{\partial u}{\partial d}(x)$ exists on $(a,b)$ and is the solution of the variational equation \textup{(\ref{var})} along $u(x)$ satisfying the boundary conditions
	\begin{align*}
		&D^{(i)}\left(x_j\right) = 0,\; 0\le i\le m_j-1,\;1 \leq j \leq k-1, \\
		&D^{(i)}\left(x_k\right) + \int_c^dpD(x)dx = pu(d),\; 0 \leq i \leq m_k - 1.
	\end{align*}
	
	\item [(e)]  $P(x):=\frac{\partial u}{\partial p}(x)$ exists on $(a,b)$ and is the solution of the variational equation \textup{(\ref{var})} along $u(x)$ satisfying the boundary conditions
	\begin{align*}
		&P^{(i)}\left(x_j\right) = 0,\; 0\le i\le m_j-1,\;1 \leq j \leq k-1, \\
		&P^{(i)}\left(x_k\right) + \int_c^dpP(x)dx = -\int_c^du(x)dx,\; 0 \leq i \leq m_k - 1.
	\end{align*}
\end{enumerate}
\end{theorem}
\bigskip

\begin{proof}

We only prove part (a) as the proofs of (b)-(e) follow similarly.  Fix integers $1 \leq l\leq k-1$ and $0\le r\le m_l-1.$ We consider $Y_{rl}(x)=\frac{\partial u}{\partial y_{rl}}(x).$  Since the argument for the case of $Y_{ik}(x)=\frac{\partial u}{\partial y_{ik}},\;\;0\le i\le m_k-1,$ is similar, we omit its proof.

To ease the burdensome notation and realizing that all boundary data are fixed except $y_{rl}$, we denote $u(x, x_1,\ldots, x_k,\\ y_{01}, \ldots, y_{rl}, \ldots, y_{m_k -1, k},p,c,d)$ by $u(x, y_{rl})$. Let $\delta > 0$ be as in Theorem \ref{contdep} with $0 \leq |h| \leq \delta$, and define the difference quotient for $y_{rl}$ by
	\begin{displaymath}
  Y_{rlh}(x) = \frac{1}{h}\left[u\left(x, y_{rl} + h\right) - u\left(x, y_{rl}\right)\right].
\end{displaymath}

First, we inspect the boundary conditions for $Y_{rlh}$. Note that for every $h \neq 0$ and $0 \leq i \leq m_j - 1,\; $ $1 \leq j \leq k-1,\; j\neq l,$
\begin{align*}
  Y_{rlh}^{(i)}\left(x_j\right) &= \frac{1}{h}\left[u^{(i)}\left(x_j, y_{rl} + h\right) - u^{(i)}\left(x_j, y_{rl}\right)\right] \\
  &= \frac{1}{h}\left[y_{ij} - y_{ij}\right] \\
  &= 0,
\end{align*}

for every $0 \leq i \leq m_l - 1,\;i \neq r$
\begin{align*}
	Y_{rlh}^{(i)}\left(x_l\right) &= \frac{1}{h}\left[u^{(i)}\left(x_l, y_{rl} + h\right) - u^{(i)}\left(x_l, y_{rl}\right)\right]\\
	&= \frac{1}{h}\left[y_{il} - y_{il}\right]\\
	&= 0,
\end{align*}

and 
\begin{align*}
	Y_{rlh}^{(r)}\left(x_l\right) &= \frac{1}{h}\left[u^{(r)}\left(x_l, y_{rl} + h\right) - u^{(r)}\left(x_l, y_{rl}\right)\right]\\
	&= \frac{1}{h}\left[y_{rl} + h - y_{rl}\right]\\
	&= 1.
\end{align*}

Finally, for every $0 \leq i \leq m_k -1$

\begin{align*}
	Y_{rlh}^{(i)}\left(x_k\right) + \int_c^d pY_{rlh}(x)dx &= \frac{1}{h}\left[u^{(i)}\left(x_k, y_{rl} + h\right) - u^{(i)}\left(x_k, y_{rl}\right)\right.\\
	&+\left.\int_c^d p\left(u(x,y_{rl}+h)- u(x,y_{rl})\right)dx\right]\\
	&= \frac{1}{h}\left[y_{ik} - y_{ik}\right] \\
	&= 0.
\end{align*}

Next, we show that $Y_{rlh}(x)$ is a solution of the variational equation. To that end, for $m_l \leq i \leq n-1$, let 
\begin{displaymath}
  \mu_i = u^{(i)}\left(x_l, y_{rl}\right)
\end{displaymath}

and 
\begin{displaymath}
 \nu_i = \nu_i(h) = y^{(i)}\left(x_l, y_{rl} + h\right) - \mu_i
\end{displaymath}

Note by Theorem \ref{contdep}, for $m_l\leq i\leq n-1, \;\nu_i = \nu_i (h) \to 0$ as $h\to
0.$  Using the notation of Theorem \ref{peano} for solutions of
initial value problems for (\ref{eq1}), viewing $u(x)$ as the solution of an initial value problem at $x_l,$ and denoting this solution as an IVP, i.e. $u(x) = y\left(x,x_l,y_{0l},\ldots,y_{m_l-1,l},\mu_{m_l},\ldots,\mu_{n-1}\right)$, we have

\begin{align*}
  Y_{rlh}(x) = \frac{1}{h}[ &y(x, x_l, y_{0l}, \ldots, y_{rl} + h, \ldots, y_{m_l -1,l}, \mu_{m_l} + \nu_{m_l}, \mu_{m_l + 1} + \nu_{m_l + 1}, \ldots, \mu_{n-1} + \nu_{n-1}) \\&- y(x, x_l, y_{0l}, \ldots, y_{rl}, \ldots, y_{m_l -1,l}, \mu_{m_l}, \mu_{m_l + 1}, \ldots, \mu_{n-1})].
\end{align*}

Next, by utilizing telescoping sums to vary only one component at a time, we have 
\begin{align*}
  Y_{rlh}(x) =& \frac{1}{h}[y(x, x_l, y_{0l},\dots,y_{rl} + h, \ldots,\mu_{m_l} + \nu_{m_l}, \mu_{m_l + 1} + \nu_{m_l + 1}, \ldots,\mu_{n-1} + \nu_{n-1})\\&- y(x, x_l, y_{0l},\ldots,y_{rl},\ldots, \mu_{m_l} + \nu_{m_l}, \mu_{m_l + 1} + \nu_{m_l + 1}, \ldots, \mu_{n-1} + \nu_{n-1})\\&+ y(x, x_l, y_{0l},\ldots,y_{rl},\ldots, \mu_{m_l} + \nu_{m_l}, \mu_{m_l + 1} + \nu_{m_l + 1}, \ldots, \mu_{n-1} + \nu_{n-1}) \\&
  	- y(x, x_l, y_{0l},\ldots,y_{rl}, \ldots, \mu_{m_l}, \mu_{m_l + 1} + \nu_{m_l + 1}, \ldots, \mu_{n-1} + \nu_{n-1}) \\ &+ - \cdots \\ &+ y(x, x_l, y_{0l},\ldots,y_{rl},\ldots, \mu_{m_l}, \mu_{m_l + 1}, \, \mu_{n-1} + \nu_{n-1}) \\ &+ y(x, x_l, y_{0l},\ldots, y_{rl},\ldots, \mu_{m_l}, \ldots, \mu_{n- 1})].
\end{align*}

By Theorem \ref{peano} and the Mean Value Theorem, we obtain
\begin{align*}
	Y_{rlh}(x) &= \alpha_r(x; y(x, x_l, y_{0l},\ldots,y_{rl} + \bar{h}, \ldots,\mu_{ml} + \nu_{ml}, \ldots , \mu_{n - 1} + \nu_{n-1}))\\ &+ \frac{\nu_{m_l}}{h}\alpha_{ml}(x;y(x; x_l, y_{0l},\ldots,y_{rl},\ldots, \mu_{m_l} + \bar{\nu}_{m_l},\mu_{m_l + 1} + \nu_{m_l + 1}, \ldots, \mu_{n-1} + \nu_{n-1}))\\
	&+ \cdots \\&+ \frac{\nu_{n-1}}{h}\alpha_{n-1}(x;y(x, x_l, y_{0l}, \ldots,\mu_{m_l}, \mu_{m_l + 1}, \ldots, \mu_{n-1} + \bar{\nu}_{n-1})),
\end{align*}

where for $0 \leq j \leq n-1,$ $\alpha_j(x; y(\cdot))$ is the solution of the variational equation (\ref{eq1}) along $y(\cdot)$ satisfying

\begin{displaymath}
  \alpha^{(i)}_j\left(x_l\right) = \delta_{ij}, \; 0 \leq i \leq n-1.
\end{displaymath}

Furthermore, $y_{rl} + \bar{h}$ is between $y_{rl}$ and $y_{rl} + h$, and for each $m_l\le i\le n-1$, $\mu_i + \bar{\nu_i}$ is between $\mu_i$ and $\mu_i + \nu_i$. Note that we use $y(\cdot)$ to simplify the notation.

Thus, to show $\ds\lim_{h \rightarrow 0} Y_{rlh}$ exists, it suffices to show, for $m_l \leq i \leq n-1,$ $\ds\lim_{h \rightarrow 0}$ $\frac{\nu_i}{h}$ exists. Recall that 
\begin{align*}
  &Y_{rlh}^{(i)}\left(x_j\right) = 0, \; 0 \leq i \leq m_j - 1,\; 1 \leq j \leq k-1, \;j \neq l,  \\
  &Y_{rlh}^{(i)}\left(x_k\right) + \int_c^d pY_{rlh}(x)dx = 0, \;0 \leq i \leq m_k - 1.  
 \end{align*}
 
 Hence, by substituting into the equations above and solving each for $\alpha_r,$  we create a system of $n-m_l$ equations with $n-m_l$ unknowns
 \begin{align*}
 	-\alpha_r^{(i)}\left(x_j; y(\cdot)\right) =  \frac{\nu_{m_l}}{h}\alpha^{(i)}_{m_l}\left(x_j; y(\cdot)\right)& + \cdots + \frac{\nu_{n-1}}{h}\alpha^{(i)}_{n-1}\left(x_j; y(\cdot)\right),\\& 0\le i\le m_j-1,\;1\le j\le k-1,\;j\ne l
 \end{align*}
 
 and 
 \begin{align*}
 	-\alpha_r^{(i)}\left(x_k; y(\cdot)\right) &- \int_c^d p\alpha_r\left(x;y(\cdot)\right)dx=  \frac{\nu_{m_l}}{h}\alpha^{(i)}_{m_l}\left(x_k; y(\cdot)\right)+ \int_c^d p\alpha_{m_l}(x;y(\cdot))dx \\&+ \cdots + \frac{\nu_{n-1}}{h}\alpha^{(i)}_{n-1}\left(x_k; y(\cdot)\right) + \int_c^d p\alpha_{n-1}(x;y(\cdot))dx,\;0\le i\le m_k-1.
 \end{align*}

In the system of equations above, we notice that $y(\cdot)$ is not always the same. Therefore, we consider the matrix along $y(x).$

 $$ M := \begin{pmatrix}
  	\alpha_{m_l}(x_1; y(x)) & \alpha_{m_l+1}(x_1;y(x)) & \cdots & \alpha_{n-1}(x_1; y(x))\\
  	\alpha^{'}_{m_l}(x_1; y(x)) & \alpha^{'}_{m_l+1}(x_1; y(x)) & \cdots & \alpha^{'}_{n-1}(x_1; y(x))\\
  	\vdots & \vdots & \ddots & \vdots\\
  	\alpha_{m_l}^{(m_1 -1)}(x_1; y(x)) & \alpha^{(m_1 -1)}_{m_l+1}(x_1; y(x)) & \cdots & \alpha^{(m_1 -1)}_{n-1}(x_1; y(x))\\
  	\vdots & \vdots & \ddots & \vdots\\
  	\alpha_{m_l}^{(m_{l-1} -1)}(x_{l-1}; y(x)) & \alpha^{(m_{l-1} -1)}_{m_l+1}(x_{l-1}; y(x)) & \cdots & \alpha^{(m_{l-1} -1)}_{n-1}(x_{l-1}; y(x))\\
  	\alpha_{m_l}(x_{l+1}; y(x)) & \alpha_{m_l + 1}(x_{l+1}; y(x)) & \cdots & \alpha_{n-1}(x_{l+1}; y(x))\\
  	\vdots & \vdots & \ddots & \vdots\\
  	
  	\alpha_{m_l}(x_k; y(x)) &  \alpha_{m_l+1}(x_k; y(x)) & \cdots & \alpha_{n-1}(x_k; y(x))\\
  	+ \int_c^d p\alpha_{m_l}(x;y(x))dx & +\int_c^d p\alpha_{m_l + 1}(x;y(x))dx & \cdots & +\int_c^d p\alpha_{n-1}(x;y(x))dx\\
  	
  	\vdots & \vdots & \ddots & \vdots\\
  	\alpha_{m_l}^{(m_k -1)}(x_k; y(x)) & \alpha_{m_l + 1}^{(m_k -1)}(x_k; y(x))  & \cdots & \alpha_{n-1}^{(m_k -1)}(x_k; y(x)) \\
  	+ \int_c^d p\alpha_{m_l} (x;y(x))dx & + \int_c^d p\alpha_{m_l + 1}(x;y(x))dx & \cdots & + \int_c^d p\alpha_{n-1}(x;y(x))dx
  	
  \end{pmatrix}$$
  
We claim that $\det(M) \neq 0$. Suppose to the contrary that $\det(M) = 0$. Then, there exists a linear combination of the column vectors with scalars $p_i\in\R,\;m_l\le i\le n-1$ such that at least one $p_i$ is nonzero

$$p_{m_l} \begin{pmatrix}
	\alpha_{m_l}(x_1; y(x))\\
	\alpha^{'}_{ml}(x_1; y(x))\\
	\vdots\\
	\alpha_{m_l}^{(m_{l-1} -1)}(x_{l-1}; y(x))\\
	\alpha_{m_l}(x_{l+1}; y(x))\\
	\vdots\\
	\alpha_{m_l}^{(m_k -1)}(x_k; y(x) \\ + \int_c^d p\alpha_{m_l}(x;y(x))dx
\end{pmatrix} + \cdots +  p_{n-1} \begin{pmatrix}
	\alpha_{n-1}(x_1; y(x))\\
	\alpha^{'}_{n-1}(x_1; y(x))\\
	\vdots\\
	\alpha_{n-1}^{(m_{l-1} -1)}(x_{l-1}; y(x))\\
	\alpha_{n-1}(x_{l+1}; y(x))\\
	\vdots\\
	\alpha_{n-1}^{(m_k -1)}(x_k; y(x)) \\ + \int_c^d p\alpha_{n-1}(x;y(x))dx
\end{pmatrix} = \begin{pmatrix}
	0 \\
	0\\
	\vdots \\
	0\\
	0\\
	\vdots\\
	0\\
\end{pmatrix}.$$

Set \begin{displaymath}
  w(x; y(x)) := p_{m_l}\alpha_{m_l}(x; y(x)) + \cdots + p_{n-1}\alpha_{n-1}(x;y(x)).
\end{displaymath}
Then by Theorem \ref{peano}, $w(x; y(x))$ is a nontrivial solution of (\ref{var}), but \begin{displaymath}
  w^{(i)}(x_j; y(x)) = 0,\; 0 \leq i \leq m_j -1, \;1 \leq j \leq k-1,\; j\neq l
\end{displaymath}
and 
\begin{displaymath}
   w^{(i)}(x_k; y(x)) + \int_c^d pw(x; y(x))dx= 0,\;0\le i\le m_k-1.
\end{displaymath}
 When coupled with hypothesis (v), we have $w(x; y(x)) \equiv 0$. Since each alpha function is not identically zero, $p_{m_l} = p_{m_{l+1}} = \cdots = p_{n-1} = 0$ which is a contradiction to the choice of $p_i$'s. Hence, $\det(M) \neq 0$ implying $M$ and, subsequently by Theorem \ref{contdep}, $M(h)$ have inverses. Here, $M(h)$ is the appropriately defined matrix from the system of equations using the correct $y(\cdot)$. Therefore, for each $m_l \leq i \leq n-1$, we can solve for $\frac{\nu_i}{h}$ by using Cramer's Rule.
 
 and suppressing the arguments of each $\alpha$:
 \begin{align*}\frac{\nu_i}{h} &=\frac{1}{M(h)} \times\\ &\begin{vmatrix}
 	\alpha_{m_l} &  \cdots & \alpha_{i-1} & -\alpha_r & \alpha_{i+1} & \cdots & \alpha_{n-1} \\
 	\vdots & \ddots & \vdots & \vdots & \vdots & \ddots & \vdots\\
 	\alpha_{m_l} + \int p\alpha_{m_l}&  \cdots & \alpha_{i-1} + \int p\alpha_{i - 1} & -\alpha_r - \int p\alpha_{r}& \alpha_{i+1} + \int p\alpha_{i + 1} & \cdots & \alpha_{n-1} + \int p\alpha_{n-1}
 \end{vmatrix}\end{align*}
 
Note as $h \rightarrow 0$, $\det(M(h)) \rightarrow \det(M),$ and so, for $1 \leq i \leq n-1 $, $\nu_i(h)/h \rightarrow \det(M_i) / \det(M) := B_i$ as $h \rightarrow 0$, where $M_i$ is the $n- m_l \times n- m_l$ matrix found by replacing the appropriate column of the matrix $M$ by 

\begin{align*}
\textup{col}\Big[&-\alpha_r(x_1;y(x)),\ldots,-\alpha_r^{(m_1-1)}(x_1;y(x)),\ldots,-\alpha_r(x_{l-1};y(x)),\ldots,-\alpha_r^{(m_{l-1}-1)}(x_{l-1};y(x)),\\&-\alpha_r(x_{l+1};y(x)),\ldots,-\alpha_r^{(m_{l+1}-1)}(x_{l+1};y(x)),
\ldots,-\alpha_r(x_k;y(x))-\int_c^d p\alpha_r(x;y(x))dx,\\&\ldots,-\alpha_r^{(m_k-1)}(x_k;y(x))-\int_c^d p\alpha_r(x;y(x))dx\Big].
\end{align*}

\par Now, let $Y_{rl}(x) = \ds\lim_{h \rightarrow 0}Y_{rlh}(x)$, and note by construction
 $$Y_{rl}(x) = \frac{\partial y}{\partial y_{rl}}(x)=\frac{\partial u}{\partial y_{rl}}(x).$$

Futhermore,
 $$Y_{rl}(x) = \lim_{h \rightarrow 0} Y_{rlh}(x) =\alpha_r\left(x;u(x)\right)+\sum_{i=m_l}^{n-1}B_i\alpha_i\left(x;u(x)\right)$$
which is a solution of the variational equation (\ref{var}) along $u(x)$. In addition,
\begin{align*}
&Y_{rl}^{(i)}\left(x_j\right) = \lim_{h \rightarrow 0}Y^{(i)}_{rlh}\left(x_j\right) = 0, \; 0 \leq i \leq m_j-1, \; 1 \leq j \leq k-1,j\ne l,\\
&Y_{rl}^{(i)}\left(x_l\right) = \lim_{h \rightarrow 0}Y_{rlh}^{(i)}\left(x_l\right) = 0, \; 0 \leq i \leq m_j-1, i\ne r,\\
&Y_{rl}^{(r)}\left(x_l\right) = \lim_{h \rightarrow 0}Y^{(i)}_{rlh}\left(x_l\right) = 1,\\
&Y_{rl}^{(i)}\left(x_k\right) + \int_c^d p Y_{rl}(x)\;dx = \lim_{h \rightarrow 0} \left[Y_{rlh}^{(i)}\left(x_k\right) + \int_c^d Y_{rlh}(x)\;dx\right] = 0,\;0 \leq i \leq m_k -1.
\end{align*}
\end{proof}

Finally, we note that similar to part (c) of Peano's theorem, the solutions found in (a)-(e) of the main result may be written as various combinations of one another due to the dimensionality of the solution space. We refer the reader to Corollary 4.1 in \cite{Lyons11} for an example.

\bibliographystyle{amsplain}
\bibliography{References}

\providecommand{\bysame}{\leavevmode\hbox to3em{\hrulefill}\thinspace}
\providecommand{\MR}{\relax\ifhmode\unskip\space\fi MR }
\providecommand{\MRhref}[2]{%
  \href{http://www.ams.org/mathscinet-getitem?mr=#1}{#2}
}
\providecommand{\href}[2]{#2}
\begin{thebibliography}{10}

\bibitem{BaxterLyonsNeugebauer16}
Lee~H. Baxter, Jeffrey~W. Lyons, and Jeffrey~T. Neugebauer,
  \emph{Differentiating solutions of a boundary value problem on a time scale},
  Bull. Aust. Math. Soc. \textbf{94} (2016), no.~1, 101--109. \MR{3539326}

\bibitem{BenchohraHamaniHendersonNtouyas07}
M.~Benchohra, S.~Hamani, J.~Henderson, S.~K. Ntouyas, and A.~Ouahab,
  \emph{Differentiation and differences for solutions of nonlocal boundary
  value problems for second order difference equations}, Int. J. Difference
  Equ. \textbf{2} (2007), no.~1, 37--47. \MR{2374098}

\bibitem{BenchohraHendersonLucaOuahab14}
Mouffak Benchohra, Johnny Henderson, Rodica Luca, and Abdelghani Ouahab,
  \emph{Boundary data smoothness for solutions of second order ordinary
  differential equations with integral boundary conditions}, Dynam. Systems
  Appl. \textbf{23} (2014), no.~2-3, 133--143. \MR{3241867}

\bibitem{Datta98}
Anjali Datta, \emph{Differences with respect to boundary points for right focal
  boundary conditions}, J. Differ. Equations Appl. \textbf{4} (1998), no.~6,
  571--578. \MR{1664373}

\bibitem{EhmeEloeHenderson93}
Jeffrey Ehme, Paul~W. Eloe, and Johnny Henderson, \emph{Differentiability with
  respect to boundary conditions and deviating argument for
  functional-differential systems}, Differential Equations Dynam. Systems
  \textbf{1} (1993), no.~1, 59--71. \MR{1385794}

\bibitem{EhmeHenderson96}
Jeffrey Ehme and Johnny Henderson, \emph{Functional boundary value problems and
  smoothness of solutions}, Nonlinear Anal. \textbf{26} (1996), no.~1,
  139--148. \MR{1354796}

\bibitem{EhmeLawrence00}
Jeffrey Ehme and Bonita~A. Lawrence, \emph{Linearized problems and continuous
  dependence for finite difference equations}, PanAmer. Math. J. \textbf{10}
  (2000), no.~2, 13--24. \MR{1754508}

\bibitem{Ehme93}
Jeffrey~A. Ehme, \emph{Differentiation of solutions of boundary value problems
  with respect to nonlinear boundary conditions}, J. Differential Equations
  \textbf{101} (1993), no.~1, 139--147. \MR{1199486}

\bibitem{EhrkeHendersonKunkelSheng07}
John Ehrke, Johnny Henderson, Curtis Kunkel, and Qin Sheng, \emph{Boundary data
  smoothness for solutions of nonlocal boundary value problems for second order
  differential equations}, J. Math. Anal. Appl. \textbf{333} (2007), no.~1,
  191--203. \MR{2323485}

\bibitem{Hartman64}
Philip Hartman, \emph{Ordinary differential equations}, John Wiley \& Sons,
  Inc., New York-London-Sydney, 1964. \MR{0171038}

\bibitem{Henderson84}
Johnny Henderson, \emph{Right focal point boundary value problems for ordinary
  differential equations and variational equations}, J. Math. Anal. Appl.
  \textbf{98} (1984), no.~2, 363--377. \MR{730513}

\bibitem{Henderson87}
\bysame, \emph{Disconjugacy, disfocality, and differentiation with respect to
  boundary conditions}, J. Math. Anal. Appl. \textbf{121} (1987), no.~1, 1--9.
  \MR{869514}

\bibitem{HendersonHopkinsKimLyons08}
Johnny Henderson, Britney Hopkins, Eugenie Kim, and Jeffrey Lyons,
  \emph{Boundary data smoothness for solutions of nonlocal boundary value
  problems for {$n$}-th order differential equations}, Involve \textbf{1}
  (2008), no.~2, 167--181. \MR{2429657}

\bibitem{HendersonHornHoward94}
Johnny Henderson, Mark Horn, and Laura Howard, \emph{Differentiation of
  solutions of difference equations with respect to boundary values and
  parameters}, Comm. Appl. Nonlinear Anal. \textbf{1} (1994), no.~2, 47--60.
  \MR{1280114}

\bibitem{HendersonJiang15}
Johnny Henderson and Xuewei Jiang, \emph{Differentiation with respect to
  parameters of solutions of nonlocal boundary value problems for difference
  equations}, Involve \textbf{8} (2015), no.~4, 629--636. \MR{3366014}

\bibitem{HendersonKarnaTisdell05}
Johnny Henderson, Basant Karna, and Christopher~C. Tisdell, \emph{Existence of
  solutions for three-point boundary value problems for second order
  equations}, Proc. Amer. Math. Soc. \textbf{133} (2005), no.~5, 1365--1369.
  \MR{2111960}

\bibitem{HendersonLee91}
Johnny Henderson and Linda Lee, \emph{Continuous dependence and differentiation
  of solutions of finite difference equations}, Internat. J. Math. Math. Sci.
  \textbf{14} (1991), no.~4, 747--756. \MR{1125427}

\bibitem{HendersonTisdell04}
Johnny Henderson and Christopher~C. Tisdell, \emph{Boundary data smoothness for
  solutions of three point boundary value problems for second order ordinary
  differential equations}, Z. Anal. Anwendungen \textbf{23} (2004), no.~3,
  631--640. \MR{2094600}

\bibitem{HopkinsKimLyonsSpeer09}
Britney Hopkins, Eugenie Kim, Jeffrey Lyons, and Kaitlin Speer, \emph{Boundary
  data smoothness for solutions of nonlocal boundary value problems for second
  order difference equations}, Comm. Appl. Nonlinear Anal. \textbf{16} (2009),
  no.~2, 1--12. \MR{2526876}

\bibitem{JansonJumanLyons14}
Alfredo~F. Janson, Bibi~T. Juman, and Jeffrey~W. Lyons, \emph{The connection
  between variational equations and solutions of second order nonlocal integral
  boundary value problems}, Dynam. Systems Appl. \textbf{23} (2014), no.~2-3,
  493--503. \MR{3241893}

\bibitem{Lawrence02}
Bonita~A. Lawrence, \emph{A variety of differentiability results for a
  multi-point boundary value problem}, J. Comput. Appl. Math. \textbf{141}
  (2002), no.~1-2, 237--248, Dynamic equations on time scales. \MR{1908841}

\bibitem{Lyons14}
J.~W. Lyons, \emph{Disconjugacy, differences and differentiation for solutions
  of non-local boundary value problems for {$n$}th order difference equations},
  J. Difference Equ. Appl. \textbf{20} (2014), no.~2, 296--311. \MR{3173548}

\bibitem{Lyons11}
Jeffrey~W. Lyons, \emph{Differentiation of solutions of nonlocal boundary value
  problems with respect to boundary data}, Electron. J. Qual. Theory Differ.
  Equ. (2011), No. 51, 11. \MR{2825136}

\bibitem{LyonsMajorSeabrook18}
Jeffrey~W. Lyons, Samantha~A. Major, and Kaitlyn~B. Seabrook, \emph{Continuous
  dependence and differentiating solutions of a second order boundary value
  problem with average value condition}, Involve \textbf{11} (2018), no.~1,
  95--102. \MR{3681350}

\bibitem{LyonsMiller15}
Jeffrey~W Lyons and Joseph~K Miller, \emph{The derivative of a solution to a
  second order parameter dependent boundary value problem with a nonlocal
  integral boundary condition}, Journal of Mathematics and Statistical Science
  \textbf{2015} (2015), no.~2, 43.

\bibitem{Peterson76}
Allan~C. Peterson, \emph{Comparison theorems and existence theorems for
  ordinary differential equations}, J. Math. Anal. Appl. \textbf{55} (1976),
  no.~3, 773--784. \MR{432977}

\bibitem{Peterson78}
\bysame, \emph{Existence-uniqueness for ordinary differential equations}, J.
  Math. Anal. Appl. \textbf{64} (1978), no.~1, 166--172. \MR{481258}

\bibitem{Spencer75}
James~D. Spencer, \emph{Relations between boundary value functions for a
  nonlinear differential equation and its variational equations}, Canad. Math.
  Bull. \textbf{18} (1975), no.~2, 269--276. \MR{399559}

\bibitem{Sukup75}
Dwight~V. Sukup, \emph{On the existence of solutions to multipoint boundary
  value problems}, Rocky Mountain J. Math. \textbf{6} (1976), no.~2, 357--375.
  \MR{409955}

\end{thebibliography}

\end{document}